\documentclass[12pt,reqno,twoside]{amsart}
\IfFileExists{euler.sty}{\IfFileExists{palatino.sty}{\usepackage{palatino}}{}}{}
\makeatletter
\@ifpackageloaded{palatino}{\usepackage[mathcal]{euler}}{}
\makeatother

\AtEndDocument{
  \vfill\eject
  \ifdim\paperheight=11in
   \message{RECOMMENDATION: If your paper is A4, insert the "a4paper" option 
   in the "documentclass" command.}
   \fi
}

\newtheorem{theorem}{Theorem}
\newtheorem{corollary}[theorem]{Corollary}
  
\def\beq#1\eeq{\begin{equation}#1\end{equation}}
\makeatletter
\@ifpackageloaded{euler}{
 \newcommand{\fchoice}[2]{#1}
 \def\B@R#1#2{\raisebox{-.07ex}{$#1#2$}\mkern-6mu}
 \renewcommand{\hbar}{{\mspace{1mu}\mathpalette\B@R{\mathchar'26}h}}
 \newcommand{\onto}{\to\mkern-14mu\to}
 \renewcommand{\cong}{\stackrel{\raise1pt\hbox{$\sim$}}{\boldsymbol{\smash=}}}
 \newcommand{\TM}{T\M}
}{
 \AtEndDocument{\IfFileExists{palatino.sty}{\message{RECOMMENDATION: Get the 
   "euler" package from tug.ctan.org.}}{\message{RECOMMENDATION: Get the 
   "palatino" and "euler" packages from tug.ctan.org. }}}
 \newcommand{\fchoice}[2]{#2}
 \renewcommand{\hbar}{{\mathchar'26\mkern-9muh}}
 \DeclareMathSymbol{\onto}{\mathrel}{AMSa}{"10}
 \newcommand{\TM}{T\!\M}
}
\makeatother
\newcommand{\A}{\mathcal{A}}
\renewcommand{\AA}{\mathbb{A}}
\newcommand{\AS}{\AA^\mathrm{S}}
\newcommand{\Ah}{\AA^{\!\hbar}}

\newcommand{\C}{\mathcal{C}}
\newcommand{\M}{\mathcal{M}}
\newcommand{\co}{\mathbb{C}}
\newcommand{\R}{\mathbb{R}}

\newcommand{\N}{\mathbb{N}}
\newcommand{\cs}{$\mathrm{C}^*$}

\newcommand{\into}{\hookrightarrow}
\newcommand{\tq}{T}
\newcommand{\gq}{Q}
\newcommand{\back}{Q^{\mathrm{inv}}}
\newcommand{\Po}{\mathcal{P}}
\newcommand{\Gh}{\Gamma_{\!\mathrm{hol}}}
 
\newcommand{\HN}{\mathcal{H}_N} 
\newcommand{\LN}{L_N} 

\DeclareMathOperator{\td}{td}
\newcommand{\Ahat}{\hat{A}}
\DeclareMathOperator{\ch}{ch}
\DeclareMathOperator{\tr}{tr}
\DeclareMathOperator{\Tr}{Tr}
\DeclareMathOperator{\der}{der}
\DeclareMathOperator{\id}{id}
\DeclareMathOperator{\End}{End}
\newcommand{\Kahler}{K\"ahler} 
\newcommand{\Toplitz}{T\"oplitz}
\newcommand{\norm}[1]{\lVert#1\rVert} 
\newcommand{\starh}{*_\hbar}
\newcommand{\W}{\mathbb{W}^\hbar}

\begin{document}
	
\begin{flushright}
\vspace*{-.4in}
\begin{tabular}{l}
\textsf{\small CGPG-98/11-1}\\
\textsf{\small math.QA/9811049}\\
\end{tabular}
\vspace{0.25in}
\end{flushright}

\title[Geometric and Deformation Quantization]{The Correspondence 
       Between\\ Geometric Quantization and\\ Formal 
       Deformation Quantization}
\author{Eli Hawkins}
\maketitle
\vspace{-4ex}
\begin{center}
\small\emph{\small Center for Gravitational Physics and Geometry}\\
\emph{\small The Pennsylvania State University,
University Park, PA 16802}\\
{\small E-mail: mrmuon@phys.psu.edu}\\
\end{center}

\begin{abstract}
Using the classification of formal deformation quantizations, and the 
formal, algebraic index theorem, I give a simple proof as to which 
formal deformation quantization (modulo isomorphism) is derived from a 
given geometric quantization.
\end{abstract}

\section{Introduction}
\label{intro}
There are two principal mathematical notions of ``quantization''\!. Both share 
as a starting point the idea (from physics) that the product of functions on a 
manifold is deformed with a parameter $\hbar$ in such a way that 
the commutator is given to leading order by the Poisson bracket as 
\beq
[f,g]_- = -i\hbar \,\{f,g\} + \mathcal{O}^2(\hbar)
\mbox.\label{commutator}\eeq

One theory, geometric quantization, gives concrete procedures 
for constructing a \cs-algebra for each (allowed) value of $\hbar$. In 
the limit as $\hbar\to0$, each of these algebras can be linearly identified 
with the ordinary algebra of continuous functions. It is in this 
approximate sense that the elements of the algebra can be thought of 
as being fixed while the product changes and satisfies Eq.~\eqref{commutator}.

In the other theory, (formal) deformation quantization (see 
\cite{bffls,wei1}), Eq.~\eqref{commutator} is taken to suggest an expansion 
in powers of $\hbar$.  The $\hbar$-dependent product is expressed as a 
power series in $\hbar$.  This power series does not, however, converge for 
most smooth functions; hence, $\hbar$ can only be taken as a formal 
parameter and cannot be given a specific, nonzero value.

Both these theories were originally intended to address the physical 
problem of quantizing the phase space of a physical system.  Physically the 
value of $\hbar$ is not variable (in fact $\hbar\approx 
10^{-27}\,\mathrm{g\,cm^2/sec}$), so deformation quantization can never be 
used to fully describe what it was originally intended to.  However, 
deformation quantization has proven fruitful as a mathematical subject.  
For instance, interesting classification results have been achieved in this 
abstract setting (see \cite{n-t2}).

As far as using deformation quantization for something like its intended 
purpose, I believe that it should be interpreted as describing the 
asymptotic behavior of a more concrete structure, such as that produced by 
geometric quantization.  Any result concerning deformation quantization 
should then have implications for this concrete version of quantization.

Essential to exploiting this is an understanding of the relationship 
between deformation quantization and geometric quantization.  In principle, 
any geometric quantization can be viewed in terms of an $\hbar$-dependent 
product which can then be asymptotically expanded to yield a deformation 
quantization.  Since this procedure is laborious at best, a shortcut to 
understanding what deformation quantization this gives is desirable.  A 
comparison of my results on ``quantization of vector bundles'' \cite{haw2} 
with an index theorem in deformation quantization theory achieves this.

\section{Geometric Quantization}
\label{geometric}
My paper \cite{haw2} is a reference for everything in this section.  Let 
$\M$ be a compact, connected, \Kahler\ manifold.  Let $L$ be a Hermitean 
line bundle with curvature equal to the symplectic form $\omega$ on 
$\M$ (which implies that $\frac\omega{2\pi}$ is integral). 
Let $L_0$ a holomorphic line bundle with an inner product on sections 
(making $\Gamma(\M,L_0)$ a pre-Hilbert space). From these we can construct 
a sequence of holomorphic line bundles $\LN:=L_0\otimes L^{\otimes N}$ 
which also have an inner product on sections.

For any $N$ the space 
$\HN:=\Gh(\M,\LN)$ of holomorphic sections of $\LN$ is a finite-dimensional 
Hilbert subspace of $L^2(\M,\LN)$. For $N$ sufficiently negative, 
$\LN$ is a negative line bundle and so $\HN=0$. Using these Hilbert 
spaces, we can define the matrix algebras $\A_N:=\End(\HN)$.

The \Toplitz\ quantization maps are completely positive maps 
$\tq_N:\C(\M)\onto\A_N$. For any continuous function $f\in\C(\M)$, the 
action of the operator $\tq_N(f)$ is defined by  multiplying a section in 
$\HN$ by $f$ and then orthogonally projecting back to $\HN$.

The geometric quantization maps $\gq_N:\C(\M)\onto\A_N$, are also 
completely positive. They can be expressed as $\gq_N(f) = 
\tq_N(f+\Delta f/2N)$, where $\Delta$ is the Laplacian.

Both of these systems of maps can be assembled into a direct-product map 
$\tq$~(respectively $\gq$)~$: \C(\M) \to \prod_{N} \A_N$, where $\prod_N 
\A_N$ is the \cs-algebraic direct product of the $\A_N$'s.  Define $\AA$ to 
be the \cs-algebra spanned by the image of $\tq$ (or $\gq$, the result is 
the same) and the \cs-algebraic direct sum $\bigoplus_N \A_N$.  This 
algebra $\AA$ is in fact the algebra of continuous sections of 
$\A_{\hat\N}$, a continuous field of \cs-algebras.  The base space of 
$\A_{\hat\N}$ is $\hat\N:=\{1,2,\ldots,\infty\}$, the one-point 
compactification of the natural numbers (or, better, the set of $N$ such 
that $\HN\neq0$).  The fiber over any finite $N\in\N$ is $\A_N$; the fiber 
over $\infty\in\hat\N$ is $\C(\M)$.

Define $\Po:\AA\onto\C(\M)$ to be the evaluation of sections at 
$\infty\in\hat\N$.  Define the partial traces $\tr_N : \AA \to \co$ by 
letting $\tr_N (a)$ be the trace of the section $a$ evaluated at $N$.  In 
\cite{haw2} I proved:
\begin{theorem}\label{GQindex}
Let $e=e^2\in M_m[\C(\M)]$ and $\tilde e=\tilde e^2\in M_m[\AA]$ be 
idempotent matrices such that $\Po(\tilde e) = e$.  For $N$ sufficiently 
large,
\beq
\tr_N \tilde e = \int_\M \!\!  \ch e \wedge \td \TM \wedge e^{c_1(L_0) + 
N\omega/2\pi}
\mbox.\label{GQiform}\eeq
Here $\ch e$ is the Chern character of the bundle determined by $e$, and 
$\td \TM$ is the Todd class of the holomorphic tangent bundle of $\M$.
\end{theorem}

\section{Formal Deformation Quantization}
\label{deformation}
Let $\M$ again be a symplectic manifold.  A (formal) deformation 
quantization of $\M$ (see \cite{wei1}) is an algebra $\Ah$ which (as a 
vector space) is identified with $\C^\infty(\M)[[\hbar]]$, the space of 
formal power series in $\hbar$ with coefficients in the smooth functions 
over $\M$.  Denote the $\Ah$-product by $\starh$ and the 
$\C^\infty(\M)[[\hbar]]$-product by apposition (e.~g., $fg$).  The product 
$\starh$ is given by a formal power series
\beq
f\starh g = fg + \sum_{k=1}^\infty (-i\hbar)^k \varphi_k(f,g)
\mbox.\label{star.series}\eeq
This is required to be associative and $\co[[\hbar]]$-linear.  It is 
required to satisfy $f\starh1=f$ and $f^*\starh g^* = (g\starh f)^*$ where 
the complex conjugate $\hbar^*=\hbar$.  The only condition involving the 
symplectic form is the restatement of Eq.~\eqref{commutator},
\[
f\starh g - g \starh f \equiv -i\hbar \,\{f,g\} \mod \hbar^2
\mbox,\]
or equivalently $\varphi_1(f,g)-\varphi_1(g,f)=\{f,g\}$, where 
$\{f,g\}$ is the Poisson bracket. Finally there is a 
(perhaps unnecessary) locality condition that each $\varphi_k$ is a 
bidifferential operator.

The archetypal example of a deformation quantization is the Moyal-Weyl 
deformation on a symplectic vector space $\M=\R^{2n}$.  Let $\mathfrak m : 
\C^\infty(\R^{2n})\otimes\C^\infty(\R^{2n})\to\C^\infty(R^{2n})$ be the 
(ordinary) multiplication map, and $\pi$ the Poisson bivector, regarded as 
a differential operator on $\C^\infty(\R^{2n})\otimes\C^\infty(\R^{2n})$ 
(so that $\mathfrak m\circ\pi(f\otimes g) = \{f,g\}$).  The Weyl product is
\begin{align}
f\starh g &:= \mathfrak m \circ \exp[-\tfrac{i\hbar}2\pi](f\otimes g)
\label{Weyl.product}\\&\:= fg - \tfrac{i\hbar}2 \{f,g\} + \dots
\nonumber\mbox.\end{align}
The formal Weyl algebra $\W$ is related to this, essentially by taking germs 
of functions about $0\in\R^{2n}$. It is constructed using formal power 
series $\co[[\R^{2n}]]$ in place of smooth functions; in other words, $\W$ 
is $\co[[\R^{2n},\hbar]]$ with the product \eqref{Weyl.product}.

Over any manifold, we can construct a bundle $\co[[\TM]]$ of formal power 
series over each fiber of the tangent bundle.  A Leibniz connection over a 
bundle of algebras, such as this, is one satisfying the Leibniz rule with 
respect to the product of sections.  The constant sections of $\co[[\TM]]$ 
with respect to a flat Leibniz connection are naturally identified with the 
smooth functions on $\M$.  The value of a section at some $x\in\M$ is a 
Taylor expansion about $x$ of the corresponding function.

Every fiber of the tangent bundle of a symplectic manifold is a symplectic 
vector space.  From this, we can construct a bundle $\W\M$ of formal Weyl 
algebras such that the fiber over $x\in\M$ is the formal Weyl algebra 
constructed on $T_x\M$.  Note that the order $\hbar^0$ part is just 
$\W\M/\hbar = \co[[\TM]]$.  The structure Lie algebra for a Leibniz 
connection on $\W\M$ is $\der\W\!$, the derivations of the typical fiber.

A flat, Leibniz connection $\nabla$ on $\W\M$ is known as a Fedosov 
connection.  The algebra $\Ah$ of $\nabla$-constant sections of $\W\M$ is a 
deformation quantization of $\C^\infty(\M)$.  Fedosov connections always 
exist \cite{fed2}, and, in fact, any deformation quantization can be 
constructed in this way \cite{n-t1}.

The space $\mathfrak g := \hbar^{-1}\W$ (series with an order $\hbar^{-1}$ 
term allowed) is a Lie algebra with the commutator as a Lie bracket.  
Indeed, $\mathfrak g$ acts on $\W$ by derivations and in fact gives all 
derivations of $\W\!$.  It is thus a central extension
\beq
0\to \hbar^{-1}\co[[\hbar]] \to \mathfrak g \to \der \W \to 0
\mbox.\label{central}\eeq

Since the Fedosov connection $\nabla$ is a $\der \W$-connection, it can be 
lifted to a $\mathfrak g$-connection $\tilde \nabla\!$ using 
Eq.~\eqref{central}.  The flatness of $\nabla$ implies that the curvature 
of $\tilde\nabla$ is central, that is 
$\tilde\nabla^2\in\hbar^{-1}\Omega^2(\M)[[\hbar]]$.  However, the lifting 
of $\nabla$ is not unique, so the curvature of $\tilde\nabla$ is not 
uniquely determined by $\nabla\!$.  Fortunately, the ambiguity is only 
modulo exact forms, so we can define $\theta := [\tilde\nabla^2]/2\pi i\in 
\hbar^{-1}H^2_\mathrm{dR}(\M)[[\hbar]]$, where brackets denote the deRham 
cohomology class.  To leading order in $\hbar$, this is given by the 
symplectic form as
\[
\theta = \frac{[\omega]}{2\pi\hbar} + \dots
\mbox.\]

The group of arbitrary automorphisms of $\Ah$ decomposes as the 
direct product of the subgroup of automorphisms preserving $\hbar$, 
with the group of formal $\hbar$ reparameterizations. 
The group of $\hbar$-preserving automorphisms is itself an extension 
of the group of symplectomorphisms ($\omega$-preserving 
diffeomorphisms) by the group of internal automorphisms. 

The cohomology class $\theta$ turns out \cite{n-t2} to classify deformation 
quantizations modulo inner automorphisms and small (connected component) 
symplectomorphisms.  The class $\theta$, modulo ``large'' 
symplectomorphisms and formal $\hbar$ reparameterization, therefore 
classifies $\Ah$ modulo isomorphisms.

For a given deformation quantization $\Ah$ of $\M$, there exists a natural 
trace (see \cite{fed3,n-t2,ros1}) $\Tr : \Ah \to \hbar^{-n}\co[[\hbar]]$.  
This is given to leading order in $\hbar$ by
\beq
\Tr f = \int_\M \frac{f\omega^n}{(2\pi\hbar)^n n!} + \dots
\mbox.\label{trace}\eeq

Using $\theta$ and this trace, a formal index theorem can be formulated 
(see \cite{n-t1} for the original, \cite{n-t3,ros1} for clarity, and also 
\cite{fed3,n-t2}).  In the case of compact $\M$, this reads:
\begin{theorem}\label{Findex}
If $e=e^2\in M_m[\C^\infty(\M)]$ and $\mathbf e=\mathbf e^2\in M_m[\Ah]$ are 
idempotents such that $\mathbf e\equiv e \mod\hbar$, then
\beq
\Tr \mathbf e = \int_\M \!\! \ch e \wedge \Ahat (\TM) \wedge e^\theta
\mbox.\label{Fiform}\eeq
\end{theorem}

\section{Comparison}
As suggested in Sec.~\ref{intro}, we can try to construct a deformation 
quantization from a geometric quantization. Let $\gq_N$ represent either the 
geometric or \Toplitz\ quantization maps. Suppose that we choose a sequence of 
maps $\back_N : \A_N\into\C^\infty(\M)$ such that $\gq_N\circ\back_N = \id$ 
and $\back_N\circ\gq_N\to\id$ as $N\to\infty$. 
Using such maps, we can pull the product on $\A_N$ back 
to $\C(\M)$ and define, 
\[
f *^N\! g := \back_N[\gq_N(f)\gq_N(g)]
\mbox.\]
If all goes well, for any smooth functions $f,g\in\C^\infty(\M)$, we can 
define $f\starh g$ by taking the asymptotic expansion of $f *^N\!  g$ as 
$N\to\infty$ and setting $N=\hbar^{-1}$.  The requirement that 
$\C^\infty(\M)[[\hbar]]$ be closed under $\starh$ is equivalent to the 
requirement that the image of each $\varphi_j$ in Eq.~\eqref{star.series} 
is in $\C^\infty(\M)$.  Associativity of $\starh$ is automatic.

The necessity of making a (somewhat arbitrary) choice of $\back_N$ in this 
construction is a rather unpleasant feature; fortunately, it can be 
eliminated.  Suppose (for any $f$, $g$, and $k$) that we take the 
difference of $f*^N\!g$ with the order $N^{-k}$ partial sum of the series 
\eqref{star.series} for $f*_{N^{-1}}g$.  The defining property of the 
asymptotic expansion is that the norm of this difference is of order 
$o^{-k}(N)$ (i.~e., the ratio with $N^{-k}$ goes to $0$ as $N\to\infty$).  
Since the quantization maps are norm-contracting, we can apply $\gq_N$ to 
the difference before taking the norm.  This cancels out the $\back_N$ and 
shows that for any $f,g\in\C^\infty(\M)$ and $k\in\N$,
\beq
\lim_{N\to\infty}
\Bigl\lVert N^k \gq_N(f)\gq_N(g) - \sum_{j=0}^{k} (-i)^j 
N^{k-j} \gq_N\left[\varphi_j(f,g)\right]\Bigr\rVert 
= 0
\mbox.\label{stardef}\eeq
Using the property that $\lim_{N\to\infty} \norm{\gq_N(f)} = \norm f$ (see 
\cite{haw2}), it is easy to verify that a $\starh$-product satisfying 
\eqref{stardef} is unique.  It is proven in \cite{gui1} that a 
$\starh$-product satisfying Eq.~\eqref{stardef} does exist and defines a 
deformation quantization, although only for the (slightly) restricted case 
of $L_0$ trivial.

Recall from Sec.~\ref{geometric} that the algebra $\AA$ is the space of 
continuous sections of $\A_{\hat\N}$, a continuous field of \cs-algebras.  
It is useful to refine $\AA$ to a space of \emph{smooth} sections.  Let 
$\AS$ be the subspace of elements in $\AA$ which have an asymptotic 
expansion as $N\to\infty$.  To be precise,
\beq
\AS := \left\{a\in\AA \Bigm| \exists f\in\C^\infty(\M)[[N^{-1}]] \; \forall 
k\in\N : \lim_{N\to\infty} N^k\norm{a - Q_N(f_{(k)})} = 0\right\}
\mbox,\label{ASdef}\eeq
where $f_{(k)}$ is the order $N^{-k}$ partial sum of $f$.  The fact that 
\eqref{stardef} can be satisfied for an algebraically closed $\starh$ 
implies that $\AS\subset\AA$ is a subalgebra.  When the $f$ in 
Eq.~\eqref{ASdef} exists, it is unique, so there is a well defined 
asymptotic expansion homomorphism $\iota : \AS\to\Ah\!, a\mapsto f$.

The subalgebra $\AS$ is holomorphically closed.  To verify this, it is 
sufficient to check that if $a\in\AS$ and $F:\co\to\co$ is holomorphic on 
the disc of radius $\norm{a}$, then $F(a)$ has the asymptotic expansion 
$F[\iota(a)]$.

Quotienting by $\hbar^k$ is a way of only dealing with the algebras $\AS$ 
and $\Ah$ to \fchoice{\linebreak}{}order $\hbar^{k-1}$.  When we do this, 
the homomorphism $\iota$ induces an isomorphism $\AS/\hbar^k\!\AS 
\widetilde\longrightarrow \Ah/\hbar^k\!$.  This leads to an invariant 
construction for $\Ah\!$.  Note that $\AS_0:=\hbar \AS$ is simply the ideal 
of sections in $\AS$ which vanish at $\hbar=0$.  Taking powers of this 
gives a nested sequence of ideals $(\AS_0)^k = \hbar^k\AS\!$.  We can 
express $\Ah$ as an algebraic inverse limit
\beq
\Ah = \varprojlim \AS\!/(\AS_0)^k
\label{invlim}\mbox,\eeq
using the obvious projections $\AS/(\AS_0)^{k+1}\onto\AS/(\AS_0)^k\!$.  The 
homomorphism $\iota$ is recovered canonically from this inverse limit 
construction.  This construction only depends on the specification of 
$\AS\subset\AA$, which is the same for both the geometric and \Toplitz\ 
quantization maps.

In light of the classification of deformation quantizations by cohomology 
classes, the obvious question now is: If a deformation quantization can be 
successfully constructed from a geometric quantization, then what is 
$\theta$?  This question is easily answered by comparing Theorems 
\ref{GQindex} and \ref{Findex}.

\begin{theorem}
	\label{main.thm}
	For the deformation quantization derived from the geometric quantization 
	of a compact, \Kahler\ manifold $\M$, the classifying cohomology class is
	\[
	\theta = \frac{[\omega]}{2\pi\hbar} + c_1 (L_0) + \tfrac12 c_1(\TM)
	\mbox.\]
\end{theorem}
\begin{proof}
	Theorem \ref{GQindex} shows that $\tr_N 1$ grows as a polynomial in $N$. 
	With the inequality $\lvert\tr_N a\rvert \leq \norm{a}_N \tr_N 1$, 
	this shows that for any $a\in \ker \iota$ (that is, $a\sim0$), 
	$\tr_N a \sim 0$. This shows that the asymptotic expansion of 
	$\tr_N$ gives a well defined, $\co[\hbar]$-linear trace on $\text{Im}~
	\iota\subset \Ah$.  Because of \eqref{invlim}, this extends uniquely to all 
	of $\Ah$ and must therefore be proportional to the trace $\Tr$ (by the 
	uniqueness of $\Tr$, see \cite{n-t2}).  To be precise, for any $a\in\AS\!$, 
	the asymptotic expansion of $\tr_N a$ as $N=\hbar^{-1}\to\infty$ is
	\[
	\tr_N a \sim \beta \Tr[\iota(a)]
	\mbox,\] 
	where $\beta \in\co[[\hbar]]$ is independent of $a$. Of course, this 
	extends to matrices over $\AS\!$.

	Choose any $a\in M_m[\AS]$ such that $\Po(a)=e$. This is idempotent 
	modulo $\AS_0$ in the sense that $a^2-a\in M_m[\AS_0]$. Since $\AS$ is 
	holomorphically closed, we can use a standard contour 
	integral trick to construct from $a$ and idempotent $\tilde e\in 
	M_m[\AS]$ such that $\tilde e - a\in M_m[\AS_0]$. Hence, 
	$\Po(\tilde e) = e$, which is the hypothesis of Thm.~\ref{GQindex}.
	Equation \eqref{GQiform} gives an exact polynomial expression for 
	$\tr_N \tilde e$ for $N$ sufficiently large; this polynomial \emph{is} 
	the asymptotic expansion of $\tr_N \tilde e$.

The idempotent $\mathbf e := \iota(\tilde e)$ satisfies the hypothesis of 
Thm.~\ref{Findex} that $\mathbf e \equiv e \mod \hbar$, 
so $\Tr[\iota(\tilde e)]$ is given by Eq.~\eqref{Fiform}. 

Combining these results gives,
\[
\beta \int_\M \!\! \ch e \wedge \Ahat (\TM) \wedge e^\theta 
= \int_\M \!\! \ch e \wedge \td \TM \wedge e^{c_1(L_0) + \omega/2\pi\hbar}
\mbox.\]
Recall that the $\Ahat$ and Todd classes are related by $\td \TM = 
e^{\frac12 c_1 (\TM)}\wedge \Ahat (\TM)$, where $c_1$ is the first 
Chern class. Now, noting that the possible values 
of $\ch(e)$ span $H^*_\mathrm{dR}(\M)$, and that $\Ahat (\TM)$ 
is invertible, this gives
\[
\theta + \ln \beta =  c_1(L_0) + \frac{[\omega]}{2\pi\hbar} + \tfrac12 c_1(\TM)
\mbox.\]
All terms of this equation are of degree $2$ except for $\ln\beta$ 
which is of degree $0$; therefore, $\ln\beta=0$.
\end{proof}

\begin{corollary}
	For any $a\in\AS\!$, 
	\[
	 \tr_N a \sim \Tr[\iota(a)]
	\]
	as $N=\hbar^{-1}\to\infty$.
\end{corollary}

\subsection*{Acknowledgments}
I wish to thank Nigel Higson and Boris Tsygan for their advice. 
This material is based upon work supported in part under a National Science 
Foundation Graduate Fellowship. Also supported in part by NSF grant 
PHY95-14240 and by the Eberly Research Fund of the Pennsylvania State 
University.


\begin{thebibliography}{10}
\newcommand{\aF}{\textsc}
\newcommand{\bF}{\textbf}
\newcommand{\tF}{\emph}

\bibitem{bffls}
\aF{Bayen, F., Flato, M., Fronsdal, C., Lichnerowicz, A., 
Sternheimer, D.}: \tF{Deformation Theory and Quantization}. Ann.\ 
Phys.\ \textbf{111} (1977), pp.~61--151. 

\bibitem{fed2}
\aF{Fedosov, B.\ V.}:
\tF{A simple geometrical construction of deformation quantization}.  
J.\ Diff.\ Geom.\ \textbf{40} (1994), no. 2, pp.~213--238. 

\bibitem{fed3}
\aF{Fedosov, B.\ V.}: \bF{Deformation quantization and index theory.} 
Math.\ Top., 9. Berlin: Akademie Verlag, 1996. 

\bibitem{gui1}
\aF{Guillemin, V.}: \tF{Star Products on Compact Pre-quantizable 
Symplectic Manifolds}. Lett.\ Math.\ Phys.\ \textbf{35} (1995), pp.~85-89.

\bibitem{haw2}
\aF{Hawkins, E.}: \tF{Geometric Quantization of Vector Bundles}. 
E-print, math.QA/9808116.

\bibitem{n-t1}
\aF{Nest, R., Tsygan, B.}: \tF{Algebraic Index Theorem}. Comm.\ 
Math.\ Phys.\ \textbf{172} (1995), 2 , pp.~223--262.

\bibitem{n-t2}
\aF{Nest, R., Tsygan, B.}: \tF{Algebraic Index Theorem for Families}. 
Adv.\ Math.\ \textbf{113} (1995), 2, pp.~151--205.

\bibitem{n-t3}
\aF{Nest, R., Tsygan, B.}: \tF{Formal Versus Analytic Index Theorems}. 
IMRN 11 (1996).

\bibitem{ros1}
\aF{Rosenberg, J. M.}: Review of \cite{n-t1,n-t2}. Math.\ Rev.\ 
96j:58163a.

\bibitem{wei1}
\aF{Weinstein, A.}: \tF{Deformation Quantization}. S\'{e}minaire Bourbaki, 
Vol.\ 1993/94. Ast\'{e}risque No.~227 (1995), Exp.\ No.~789, 5, 389--409.

\end{thebibliography}
\end{document}